\documentclass[11pt]{amsart}
\usepackage{amsmath,amssymb,latexsym}
\usepackage[T1]{fontenc}
\usepackage{enumerate}

\def\F{\mathfrak F}

\def\rk{\mathrm{rk}}

\def\Fix{\mathrm{Fix}}

\def\Ind#1#2{#1\setbox0=\hbox{$#1x$}\kern\wd0\hbox to 0pt{\hss$#1\mid$\hss}
\lower.9\ht0\hbox to 0pt{\hss$#1\smile$\hss}\kern\wd0}

\def\Notind#1#2{#1\setbox0=\hbox{$#1x$}\kern\wd0\hbox to 0pt{\mathchardef
\nn="3236\hss$#1\nn$\kern1.4\wd0\hss}\hbox to 0pt{\hss$#1\mid$\hss}\lower.9\ht0
\hbox to 0pt{\hss$#1\smile$\hss}\kern\wd0}

\theoremstyle{plain}
\newtheorem{theorem}{Theorem}[section]
\newtheorem{prop}[theorem]{Proposition}
\newtheorem{fact}[theorem]{Fact}
\newtheorem{lemma}[theorem]{Lemma}
\newtheorem{cor}[theorem]{Corollary}

\theoremstyle{definition}
\newtheorem{defn}[theorem]{Definition}
\newtheorem{remark}[theorem]{Remark}

\def\pf{\par\noindent{\em Proof. }}

\title{Division rings with ranks}
\author{Nadja Hempel and Daniel Palac\'in}\thanks{The first author was supported by ValCoMo
(ANR-13-BS01-0006). The second author was supported by the projects SFB 878 and MTM2014-59178-P}
\address{Institut Camille Jordan; Universit\'e Lyon 1; 43 bd du 11 novembre 1918, 69622 Villeurbanne Cedex, France}
\email{hempel@math.univ-lyon1.fr}

\address{Mathematisches Institut, Universitat M\" unster, Einsteinstrasse 62,
48149 M\"unster, Germany}

\email{daniel.palacin@uni-muenster.de}

\begin{document}

\begin{abstract}
Any superrosy division ring (i.e. a division ring equipped with an abstract notion of rank) is shown to be centrally finite. Furthermore, division rings satisfying a generalized chain condition on definable subgroups are studied. In particular, a division ring of burden $n$ has dimension at most $n$ over its center, and any definable group of definable automorphisms of a field of burden $n$ has size at most $n$. Additionally, interpretable division rings in o-minimal structures are shown to be algebraically closed, real closed or the quaternions over a real closed field.
\end{abstract}

\maketitle

\section{Introduction}

An important aspect in model theory is to analyze algebraic properties of mathematical objects, such as groups and rings, which satisfy certain combinatorial behavior on their definable sets. Some of these combinatorial patterns (e.g. $\omega$-stability, stability, simplicity) yield the existence of suitable chain conditions among definable groups as well as well-behaved rank functions among definable sets. These are important tools to study algebraic properties of groups and rings.

A milestone in classifying fields from a model-theoretic point of view is a result of Macintyre \cite{Mac} which states that any infinite $\omega$-stable field is algebraically closed. This was generalized to the superstable context by Cherlin and Shelah \cite{CherShe} and therefore, by previous work of Cherlin \cite{Cherlin}, in fact any infinite superstable division ring is an algebraically closed field. Moreover, Pillay, Scanlon and Wagner \cite{PSW} showed that a wider class of division rings, namely the supersimple ones, are commutative and have indeed trivial Brauer group.
In all these cases the existence of a suitable ordinal-valued rank function plays an essential role. A more general context in which a similar rank function is present is the class of superrosy division rings which includes, besides all before mentioned fields, division rings interpretable in o-minimal structures. In the definable context, such division rings were first analyzed by Otero, Peterzil and Pillay \cite{OPP} and later by Peterzil and Steinhorn \cite{PS}. They were characterized to be either algebraically closed fields, real closed fields or the  quaternions over a real closed field. In particular, there exist non commutative superrosy division rings.

The first part of this paper is dedicated to the study of superrosy division rings using the aforementioned rank function from a purely axiomatic point of view. We show that any infinite superrosy division ring is centrally finite (Theorem \ref{ThmMain}). This gives a uniform treatment to all previous mentioned cases. 

One consequence of the presence of a well-behaved rank function is the non-existence of an infinite descending chain of definable subgroups, each having infinite index in its predecessor. Weakenings of such a forbidden pattern appear naturally in wider model theoretic classes, such as simple or NTP$_2$ theories. In these frameworks, a more general notion of rank, namely {\em weight} in simple theories and {\em burden} in NTP$_2$ theories, can be defined. Our aim is to analyze division rings in these contexts. For instance, using machinery from simplicity theory and classical results on division rings such as the Cartan-Brauer-Hua theorem, we show that a division ring with a simple theory of weight one must be commutative (Corollary \ref{CorWeight1}). Furthermore, we study division rings of finite burden, using only a generalized chain condition on definable subgroups. These include division rings interpretable in o-minimal structures and in simple structures of finite rank (e.g. any pseudofinite field or more generally any perfect PAC field with small Galois group). Other examples, which are neither o-minimal nor simple,  are pseudo real closed fields with small Galois group \cite{Sam}. Moreover, a recent result of Chernikov and Simon \cite{CS} gives rise to numerous examples: any ultra-product of the $p$-adics,  $\mathbb R((t))$, $\mathbb R((\mathbb Q))$, $\mathbb C((t))$ and $F((t^{\mathbb Q}))$ where $F$ is any perfect PAC field with small Galois group, among others. Indeed, all these examples have burden one.  We show that any division ring of burden $n$ has dimension at most $n$ over its center, and in particular a division ring of burden one is commutative. Furthermore, we prove that a field of burden $n$ and characteristic zero has at most $n$ many definable automorphisms (Proposition \ref{PropDefAut}).

The last section includes a remark on  division rings interpretable in o-minimal structures. We apply the above results to show that any such division ring is either  algebraically closed fields, real closed fields or the  quaternions over a real closed field. This generalizes the aforementioned result of Peterzil and Steinhorn to the interpretable framework as well as gives an alternative proof.

\section{Ranks \`a la Lascar}

In this section we aim to study superrosy division rings. As being centrally finite and superrosy are both properties of the theory of a division ring, we may work (if necessary) in an saturated elementary extension. In fact, all results presented here extend to arbitrary  $\aleph_0$-saturated division rings
 carrying an ordinal-valued rank function among the definable sets (in the pure language of rings) in the imaginary expansion, i.e.
$$
\rk : \{\mbox{Definable sets}\} \rightarrow \operatorname{Ord},
$$ which satisfies the following properties:
\begin{enumerate}
\item A definable set has rank zero if and only if it is finite.
\item The rank is preserved under definable bijections.
\item The Lascar inequalities: For a definable subgroup $H$ of a definable group $G$ we have that
$$
\rk(H)+\rk(G/H)\le \rk(G)\le \rk(H) \oplus \rk(G/H),
$$
\end{enumerate}
where the function $\oplus$ is the smallest symmetric strictly increasing function $f$ among pairs of ordinals such that $f(\alpha,\beta+1)=f(\alpha,\beta)+1$. More precisely,  every ordinal $\alpha$ can be written in the Cantor normal form as a finite sum $\omega^{\alpha_1}\cdot n_1 + \ldots + \omega^{\alpha_k}\cdot n_k$ for ordinals $\alpha_1>\ldots>\alpha_k$ and natural numbers $n_1,\ldots,n_k$.  If additionally $\beta=\omega^{\alpha_1}\cdot m_1 + \ldots + \omega^{\alpha_k}\cdot m_k$, then $\alpha\oplus\beta$ is defined to be $\omega^{\alpha_1}\cdot (n_1+m_1) + \ldots + \omega^{\alpha_k}\cdot (n_k+m_k)$.

The existence of such a rank yields some immediate consequences on definable groups of a superrosy division ring, which we state below.

\begin{remark}\label{RemInj}
 Let $G$ and $H$ be two definable groups and let $f:H\rightarrow G$ be a definable group morphism. Then
$$\rk({\rm Ker} f)+ \rk({\rm Im} f) \le \rk(H) \le \rk({\rm Ker} f) \oplus \rk({\rm Im} f).
$$
Thus, if $f$ is injective, then $H$ and $G$ have the same rank if and only if ${\rm Im} f$ has finite index in $G$. In particular, if $H$ is a subgroup of $G$, then $H$ and $G$ have the same rank if and only if $H$ has finite index in $G$.
\end{remark}

\begin{remark}\label{RemUCC}
 As there is no infinite strictly descending chain of ordinals, any infinite descending chain of definable groups, each of them having infinite index in its predecessor, stabilizes after finitely many steps.
\end{remark}
\begin{defn}
Let $X$ be a definable set of rank $\omega^\alpha\cdot n+\beta$ with $\beta<\omega^\alpha$ and $n$ a positive natural number. A definable subset $Y$ of $X$ is  {\em wide} in $X$ if it has rank at least $\omega^\alpha\cdot n$ and {\em negligible} with respect to $X$ if its rank is strictly smaller than $\omega^\alpha$. If there is no confusion we simply say that $Y$ is wide or respectively negligible.
\end{defn}

\begin{lemma}\label{LemmaNeg}
Any superrosy division ring has finite dimension (as a vector space) over any definable non-negligible subdivision ring.
\end{lemma}
\pf Let $D$ be a superrosy division ring of rank $\omega^\alpha\cdot n+\beta$ with $\beta<\omega^\alpha$ and $n$ a positive natural number. Suppose, towards a contradiction, that there is a definable subdivision ring $D_0$ of rank greater or equal than $\omega^\alpha$ such that $D$ has infinite dimension over $D_0$. Thus for $1=\lambda_0,\lambda_1, \dots ,\lambda_n$ linearly independent over $D_0$ we obtain
\begin{eqnarray*}
\rk(D) & \ge & \rk\Big(\bigoplus_{i\le n} D_0\lambda_i \Big)\  \geq\ \rk (D_0\lambda_n)+ \rk \Big(\bigoplus_{i\le n} D_0\lambda_i \Big/  D_0\lambda_n \Big)    \\
& \ge & \omega^\alpha +  \rk \Big(\bigoplus_{i< n}  D_0 \lambda_i \Big)\  \geq\  \cdots \ \ge\ \sum_{i\le n} \omega^ \alpha \cdot (n+1) \\ & >& \rk(D),
\end{eqnarray*}
which yields a contradiction.
\qed

Recall that two groups $H$ and $N$ are said to be {\em commensurable} if their intersection $H\cap N$ has finite index in each of them. For the following lemma we use Schlichting's theorem \cite{sch}, generalized by Bergmann and Lenstra
\cite{BL} to obtain an invariant subgroup, up to commensurability. See also \cite[Theorem 4.2.4]{wagner}.

\begin{fact}[Schlichting's theorem]\label{sch} Let $\mathfrak F$ be
a family of uniformly commensurable subgroups of a group $G$, i.e.\ there is a natural number $k$ such that for any  $H$ and $H^*$ in $\mathfrak F$ the index $[H:H\cap H^*]$ is at most $k$. Then there
is a subgroup $N$ commensurable with any $H$ in $\mathfrak F$, which is invariant under any
automorphism of $G$ stabilizing $\mathfrak F$ setwise. 
\end{fact}

\begin{lemma}\label{LemmaGrupRank}
Any wide definable additive subgroup of a superrosy division ring has finite index.
\end{lemma}
\pf
Let $D$ be a superrosy division ring, which we may assume to be $\aleph_0$-saturated, of rank $\omega^\alpha\cdot n+\beta$ with $\beta<\omega^\alpha$ and $n$ a positive natural number. Suppose towards a contradiction that there is a definable additive subgroup $A$ of $D$ of infinite index and of rank greater or equal to $\omega^\alpha \cdot n$. Consider the family of $D^\times$-translates of $A$. As there is no infinite descending chain of definable subgroups, each of them having infinite index in its predecessor, there is a finite intersection $H$ of $D^\times$-translates of $A$ such that for any $\lambda$ in $D^\times$ the index $[H:H\cap \lambda H]$ is finite. As left translation is an automorphism,  any two translates of $H$ are commensurable and hence, a compactness argument yields that  the family $\{\lambda H: \lambda \in D^\times\}$ is uniformly commensurable. Thus, by Schlichting's theorem, we can find a $D^\times$-invariant additive subgroup $I$, i.e. an ideal of $D$, which is commensurable with $H$. As $H$ has infinite index in $D$, the ideal $I$ has to be trivial and hence $H$ is finite. In particular, it has rank zero and therefore, there exists a finite intersection $N$ of $D^\times$-translates of $A$ such that $\rk(N)\ge \omega^\alpha\cdot n$ and for any $\lambda$ in $D^\times$ either $N$ and $\lambda N$ are commensurable or $\rk(N\cap \lambda N)<\omega^\alpha\cdot n$. As $N$ is clearly not finite, there is some $\lambda$ in $D^\times$ such that $N$ and $\lambda N$ are not commensurable and hence $\rk(N\cap \lambda N)<\omega^\alpha\cdot n$. Together with the following inequality
$$
\omega^\alpha\cdot n \le \rk(N) \le \rk(N\cap \lambda N) \oplus \rk(N/N\cap \lambda N),
$$
we get that $\rk(N/N\cap \lambda N)\ge \omega^\alpha$. Hence,
\begin{align*}
\rk(N+\lambda N) & \ge \rk(\lambda N)+ \rk(N+\lambda N/\lambda N) \\
& = \rk(N) + \rk(N/N\cap \lambda N) \\
& \ge \omega^\alpha \cdot n + \omega^\alpha \\
& > \rk(D),
\end{align*}
which yields a contradiction and finishes the proof. \qed

\begin{cor}\label{CorFinInd}
Let $D$ be a superrosy division ring. If a definable group morphism from $D^+$ or $D^\times$ to $D^+$ has a negligible kernel, its image has finite index in $D^+$.
\end{cor}
\pf As the kernel is negligible, by Remark \ref{RemInj} the image is wide and thus the previous lemma yields the statement. \qed

Before stating and proving the main result of this section, we show the following lemma which holds under the mere assumption that $D^\times$ satisfies the descending chain condition on centralizers.

\begin{lemma}\label{Lemma_ZDiv}
Let $D$ be a non commutative division ring satisfying the descending chain condition on centralizers, and suppose that for any non-central element $a$ in $D$ the center $Z(D)$ is contained in $a^D-a$. Then, for any natural number $m$, any element of $Z(D)$ has an $m$-root in $D$.
\end{lemma}
\pf By Kaplansky's theorem \cite[Theorem 15.15]{lam}, we can find an element $a$ in $D$ for which none of
its powers belong to $Z(D)$. As $D$ satisfies the chain condition on
centralizers, after replacing $a$ by one of its powers, we may assume that
$C_D(a) = C_D(a^n)$ for any natural number $n$.

Suppose, towards a contradiction, that there exists a natural number $m$ and an element $c$ in $Z(D)$ with no $m$-root in $D$. As $a^m$ is non-central, by assumption there is some $x$ in $D$ such that $(a^m)^x- a^m =c$. Observe that $a$ and $a^x$ commute since
$$ C_D(a^x) = C_D(a)^x = C_D(a^m)^x = C_D((a^m)^x) = C_D(a^m) = C_D(a)$$
and so
$$(a^x a^{-1})^m - c a^{-m} = ((a^m)^x  - c) a^{-m}= 1.$$
Furthermore, as $c a^{-m} $ is also non-central, one can find an element $y$ in $D$ with $(c a^{-m} )^y - c a^{-m}  =1$. Thus,  the $m$th power $(a^x a^{-1})^m$ equals to $(c a^{-m} )^y$. As $c$ was assumed to have no $m$-root in $D$, the latter yields a contradiction. \qed

\begin{theorem}\label{ThmMain}
A division ring with a superrosy theory has finite dimension over its center.
\end{theorem}

\pf Suppose, towards a contradiction, that there is a superrosy division ring which has infinite dimension over its center and let $D$ be such a division ring of minimal rank. As $D$ is clearly infinite its rank must be non-zero.

The proof consists of a series of steps.

\noindent{\bf Step 1.} 
Any proper centralizer has finite dimension over its center:

\noindent Since any proper centralizer is a subdivision ring and thus as an additive group has infinite index in $D$, by Lascar inequalities its rank is strictly smaller than the rank of $D$. Hence, we conclude by the choice of $D$.

\noindent{\bf Step 2.} $D$ has infinite dimension over any proper centralizer. In particular, the centralizer of a non-central element is negligible:

\noindent Otherwise, it has finite dimension over some subfield by Step 1 and thus it would be finite dimensional over its center by \cite[Theorem 15.8]{lam}. The second part is an immediate consequence by Lemma \ref{LemmaNeg}.

Now, for an arbitrary element $a$, let $\delta_a:D\rightarrow D$ be the map defined as $u\mapsto \delta_a(u)=au-ua$. To easer notation, we also write $[a,u]$ for $au-ua$.

\noindent{\bf Step 3.} For a non-central element $a$, the map $\delta_a$ is a surjective derivation whose kernel is the centralizer of $a$:

\noindent Easy computations yield that it is a derivation and that its kernel is $C_D(a)$.
As the $C_D(a)$ is negligible by Step 2, the image of $\delta_a$ has finite index in $D$ by Corollary \ref{CorFinInd}. Moreover, since ${\rm Im} \, \delta_a$ is a vector space over $C_D(a)$, which is an infinite division ring by \cite[Theorem 13.10]{lam}, it is indeed equal to $D$.

\noindent{\bf Step 4.} Any non-central element is transcendental over $Z(D)$:

\noindent Given a non-central element $a$, by Step 3 we can find a non-central element $u$ such that $\delta_a(u)=-1$, i.e. $[a,u]=-1$. Additionally, as the map $\delta_u$ is a derivation whose kernel clearly contains the center, the element $a$ must be transcendental over the center since $\delta_u(a)=-\delta_a(u)=1$.  

\noindent{\bf Step 5.} The characteristic is zero:

\noindent Suppose that $D$ has characteristic $p>0$. Up to replacing a non-central element $a$ by one of its powers (which is non-central by Step 4), we may assume by the chain condition on centralizers that for any natural number $n$ we have that $C_D(a)$ is equal to $C_D(a^n)$. By Step 3, we can find an element $x$ such that $\delta_a(x)=ax-xa=-1$. Thus,  for any $i\ge 1$ conjugation by $a^ i$ yields $x^{a^i}a-x^{a^{i-1}}a=1$. Hence
$$
(ax)^{a^{p}} - ax = x^{a^{p-1}}a - ax = \sum_{i=1}^{p-1} \big( x^{a^{i}}a-x^{a^i-1}a \big) - (ax-xa) = (p-1) + 1 =0.
$$
Whence $a^p$ commutes with $ax$ and so with $x$. Therefore $x$ belongs to $C_D(a^p)$ which equals to $C_D(a)$, yielding a contradiction.

\noindent{\bf Step 6.} For any non-central element $a$ and any $x$ such that $\delta_a(x)=1$, the intersection $Z(C_D(x))\cap C_D(a)$ is equal to $Z(D)$:

\noindent Let $y$ be an element of $Z(C_D(x))\cap C_D(a)$. Thus $Z(C_D(y))$ is contained in $C_D(a)$ and thus in the kernel of $\delta_a$. As additionally  $\delta_a(x)=1$, the element $x$ must be transcendental over $Z(C_D(y))$. However, the choice of $y$ yields that $x\in C_D(y)$ and so $C_D(y)$ has infinite dimension over its center. Therefore, we obtain that $y$ belongs to $Z(D)$ by Step 1.

\noindent{\bf Step 7.} For any non-central element $a$ and any $x$ such that $\delta_a(x)=1$, the map $\delta_a$ restricted to $Z(C_D(x))$ is a derivation whose image is contained in $Z(C_D(x))$ and whose kernel is $Z(D)$:

\noindent Note first that the kernel is precisely $C_D(a)\cap Z(C_D(x))$ which is equal to the center by Step 6. To prove that $\delta_a$ restricted to $Z(C_D(x))$ induces a map to  $Z(C_D(x))$, we first see that the image of $C_D(x)$ via $\delta_a$ is contained in $C_D(x)$. To do so, let $u$ be an element of $C_D(x)$. By the Jacobi identity, we have
$$
\big[a,[u,x] \big] + \big[u,[x,a] \big]  + \big[x,[a,u] \big]  =0.
$$
As $u$ commutes with $x$, the first summand $[a,[u,x]]$ is equal to $0$ and since $[x,a]=-1$, we obtain that $[u,[x,a]]$ is $0$ as well. Therefore $[x,[a,u]]$ must be $0$ and hence $\delta_a(u)=[a,u]$ belongs to $C_D(x)$. Second, let $v$ be an arbitrary element of $Z(C_D(x))$. Note first that $\delta_a(v)$ belongs to $C_D(x)$. Furthermore, for  $u$ again in $C_D(x)$, the identity
$$
\big[u,[a,v] \big] + \big[a,[v,u] \big]  + \big[v,[u,a] \big]  =0
$$
yields, similarly as above, that $\delta_a(v)=[a,v]$ commutes with $u$ and thus, as $u$ was taken to be arbitrary in $C_D(x)$, we obtain that $\delta_a(v)$ belongs to $C_D(C_D(x))$ and thus to $Z(C_D(x))$.


\noindent{\bf Step 8.} For any non-central element $a$ and any $x$ such that $\delta_a(x)=1$, we obtain
$
Z(C_D(x))=a^{Z(C_D(x))} - a:
$

\noindent Let $\sigma_a:Z(C_D(x))^\times\rightarrow D$ be the function that maps an element $y$ to $a^y-a$. We need to show that $\sigma_a$ is a surjective group morphism from $Z(C_D(x))^\times$ onto $Z(C_D(x))^+$. To show that the image of $\sigma_a$ is contained in $Z(C_D(x))^+$, note that for any $y$ in $Z(C_D(x))^\times$ we have that
$$
\sigma_a(y)=a^y-a=y^{-1}(ay-ya)=y^{-1}\delta_a(y)
$$
which belongs to $Z(C_D(x))$  by Step 7. Moreover, as $\sigma_a(1)=0$ and for any $u$ and $v$ in $Z(C_D(x))$
\begin{eqnarray*}
\sigma_a(uv) & = & a^{uv}-a \\
& = & a^{uv} - a^v + a^v -a \\
& = & (a^u-a)^v + (a^v - a) \\
& = & (a^u-a) + (a^v - a) \\
& = & \sigma_a(u) + \sigma_a(v),
\end{eqnarray*}
the map $\sigma_a$ is a group morphism. To conclude this step, it remains to show that $\sigma_a$ is surjective. For this, note first that the kernel of $\sigma_a$ is clearly $Z(C_D(x))^\times\cap C_D(a)^\times$ which is equal to $Z(D)^\times$ by Step 6. As $x$ is transcendental over $Z(D)$ by Step 4, we get that $Z(C_D(x))$ is infinite dimensional over $Z(D)$ and thus the kernel of $\sigma_a$ is negligible with respect to $Z(C_D(x))$ by Lemma \ref{LemmaNeg}. Thus, by Corollary \ref{CorFinInd}, the image of $\sigma_a$ has finite index in $Z(C_D(x))^+$.  Since $Z(C_D(x))^+$ is divisible as the characteristic is zero by Step 5, both are indeed equal.

\noindent{\bf Step 9.} The multiplicative group $Z(D)^\times$ is divisible:

\noindent By Step 8, we deduce that for any non-central element $a$ in $D$, the center  $Z(D)$ is contained in $a^D - a$ and hence, by Lemma \ref{Lemma_ZDiv}, any central element has an $m$-root in $D$ for any natural number $m$. Therefore, as any non-central element is transcendental over $Z(D)$ by Step 4, any root of a central element must belong to the center.

\noindent{\bf Step 10.} There exists a non-central element $b$ such that for any $y$ in $Z(C_D(b))$ we have that $Z(C_D(b))=Z(C_D(y))$:

\noindent Otherwise we can find an infinite descending chain of fields
$$
Z(C_D(a_0)) \gneq Z(C_D(a_1)) \gneq \ldots \gneq Z(C_D(a_i)) \gneq \ldots
$$
with $a_{i+1}$ being a non-central element belonging to $Z(C_D(a_i))$, which yields a contradiction.

\noindent{\bf Step 11.} Final contradiction:

\noindent For this final step, let $b$ be an element given by Step 10. By Step 3 the derivation $\delta_b$ on $D$ is surjective and so we can find an element $y$ such that $\delta_b(y)=-1$. In particular $\delta_y(b)=-\delta_b(y)=1$. Hence, by Step 7 the restriction of $\delta_y$ to $Z(C_D(b))$ is surjective onto $Z(C_D(b))$. As $-b^{-1}$ belongs trivially to $Z(C_D(b))$, there is some $c$ in $Z(C_D(b))$ such that $\delta_y(c)=-b^{-1}$, i.e. $\delta_c(y)=b^{-1}$. Hence, since $c$ commutes with $b$ we obtain that
$$
1=cyb-ycb=c(yb)-(yb)c=\delta_c(yb).
$$
Thus, Step 8 yields:
$$
Z(C_D(yb))=c^{Z(C_D(yb))}-c
$$
and in particular, for any $z$ in $Z(D)$ we can find an element $t$ in $Z(C_D(yb))$ such that $z=c^t-c$.
Now, we define the map $\tau:Z(D)^+\rightarrow D^\times$ as follows:
$$
\tau(z)=b^{-1}b^t, \ \mbox{ where $t\in Z(C_D(yb))$ such that $z=c^t-c$}.
$$
Next, we see that $\tau$ is a group isomorphism from $Z(D)^+$ onto $Z(D)^\times$. First, we check that it is well-defined.
To do so, let $z$ be in $Z(D)$ and consider $t$ and $s$ in $Z(C_D(yb))$ such that $z=c^t-c$ and $z=c^s-c$. Thus $t^{-1}s$ commutes with $c$. As $Z(C_D(c))=Z(C_D(b))$ by the choice of $b$, the product $t^{-1}s$ is an element of $C_D(b)$. Therefore, the value of $\tau(z)$ does not depend on the choice of $t$. Now, we show that ${\rm Im} \, \tau$ is contained in $Z(D)^\times$. Hence, consider an arbitrary element $z$ of $Z(D)$ and let $t$ be in $Z(C_D(yb))$ such that $c^t-c=z$. Thus $t$ normalizes $Z(C_D(c))$ and so $Z(C_D(b))$. Therefore $\tau(z)=b^{-1}b^t$ belongs to $Z(C_D(b))$. Moreover, as $t$ commutes with $yb$ and so with $by$ since $yb-by=1$, we have that
\begin{eqnarray*}
[yb,\tau(z)] & = & yb^t-b^{-1}b^tyb \\
& = & b^{-1}byb^t-b^{-1}b^tyb \\
& =& b^{-1}(byb)^t-b^{-1}(byb)^t \\
& = &0.
\end{eqnarray*}
Thus $\tau(z)$ commutes with $yb$ and hence with $y$ since $\tau(z)$ belongs to $Z(C_D(b))$. Therefore, we have shown that $\tau(z)$ is in $Z(C_D(b))\cap C_D(y)$ and so it is central by Step 6. \newline
Now, to prove that $\tau$ is a group morphism, note first that $\tau(0)=1$. Moreover, consider $z$ and $z'$ in $Z(D)$ and let $t$ and $s$ be the corresponding elements from $Z(C_D(yb))$ such that $\tau(z)=b^{-1}b^t$ and $\tau(z')=b^{-1}b^s$. The choice of $t$ and $s$ yields
$$
c^{ts}-c=(c^t-c)^s+(c^s-c)=z+z'.
$$
Thus

\begin{eqnarray*}
\tau(z+z') & = & b^{-1}(ts)^{-1}bts \\
& = & (b^{-1}s^{-1}b)(b^{-1}t^{-1}bt)s \\
& =& (b^{-1}t^{-1}bt)(b^{-1}s^{-1}b)s \\
& = & \tau(z)\tau(z').
\end{eqnarray*}
Hence, it remains to show that $\tau$ is an isomorphism. To check that the kernel of $\tau$ is trivial, consider an arbitrary element $z$ of the kernel and let $t$ be in $Z(C_D(yb))$ such that $z=c^t-c$. The choice of $z$ yields that $b^{-1}b^t=1$ and so $t$ commutes with $b$. As $Z(C_D(b))=Z(C_D(c))$, we have that $t$ commutes with $c$ and thus $z=0$.
Now, by Remark \ref{RemInj} the image of $\tau$ has finite index in $Z(D)^\times$ and so it is equal to $Z(D)^\times$ since the latter is divisible by Step 9. \newline Therefore $\tau$ is an isomorphism between $Z(D)^+$ and $Z(D)^\times$ and consequently, we can find a central element $z$ such that $\tau(z)=-1$. Hence, for some $t$ in $Z(C_D(yb))$ we have that $b^t=-b$, whence
$$b^2=(-b)^2=(b^t)^2=(b^2)^t$$
and so $t$ and $b^2$ commute. On the other hand, the choice of $b$ implies that $Z(C_D(b))=Z(C_D(b^2))$, thus $t$ commutes with $b$ and so
$$
-1=\tau(z)=b^{-1}b^t=b^{-1}b=1,
$$
a blatant contradiction since the characteristic of $D$ is zero by Step 5. This final contradiction finishes the proof. \qed

\section{Chain conditions on uniformly definable subgroups}

In this section we study rosy division rings. Its unique property which we use throughout the section is the following chain condition on uniformly definable subgroups: \newline
Let $\mathfrak F$ be a family of uniformly definable subgroups, then there is a natural number $d$ such that any descending chain of intersections
$$
F_0 \geq F_0 \cap F_1 \ge \cdots \ge  F_0 \cap\ldots \cap F_i \ge \cdots
$$
with $F_i$ in $\mathfrak F$ for all $i$,  each having index at least $d$ in its predecessor, has finite length.

\begin{remark}
As the collection of centralizers of one element is a family of uniformly definable additive subgroups, and moreover division rings, the above chain condition yields the ordinary chain condition on centralizers. Namely, any infinite descending chain of centralizers stabilizes after finitely many steps.
\end{remark}

The following result is an easy adaptation of \cite[Theorem 3.5]{Cedric}.

\begin{prop}\label{PropCed}
A rosy division ring of positive characteristic has finite dimension over its center.
\end{prop}
\pf Let $D$ be a rosy division ring of positive characteristic, say $p$. Assume, as we may, that $D$ is infinite and non-commutative. By the ordinary chain condition on centralizers, we may inductively suppose that any centralizer of a non-central element has finite dimension over its center.

By Kaplansky's theorem \cite[Theorem 15.15]{lam}, we can find a non-central element $a$ of $D$ for which no power belongs to the center. Note that the centralizer of $a$ is infinite. Additionally, by the chain condition on centralizers, after replacing $a$ by one of its powers, we may assume that $C_D(a)=C_D(a^p)$. Let $\delta_a$ be the derivation of $D$ given by $\delta_a(x)=ax-xa$, which is clearly a definable $C_D(a)$-linear morphism. Since $D$ has characteristic $p$ one can easily see that
$$
\delta_{a^{p^n}}(D)=\delta_a \, \circ\stackrel{p^n}{\dots}\circ \, \delta_a(D).
$$
As any descending chain of uniformly definable additive subgroups stabilizes after finitely many steps, up to finite index, there exists a natural number $n$ for which the group $\delta_{a^{p^{n+1}}}(D)$ has finite index in $\delta_{a^{p^n}}(D)$. Thus, they are equal as both are vector spaces over the infinite division ring $C_D(a)$. Replacing $a$ by $a^{p^n}$ if necessary, we may assume that $\delta_a(D)$ and $\delta_{a^p}(D)$ coincide, and so do $\delta_a(D)$ and $(\delta_a\circ\delta_a)(D)$. Thus, given an element $x$ of $D$ there is some $y$ such that $\delta_a(x)=\delta_a(\delta_a(y))$ and so $\delta_a(y)-x$ belongs to $C_D(a)$. Hence $D$ is the sum of the two $C_D(a)$-vector spaces $\delta_a(D)$ and  $C_D(a)$. Whence, as $C_D(a)=C_D(a^p)$, the latter sum is a direct sum, i.e.
$$
D=C_D(a)\oplus \delta_a(D).
$$

Now, set $H$ to be the $C_D(a)$-vector space $\delta_a(D)$, and consider the family $\mathfrak F$ of $D^\times$-translates of $H$. By the chain condition applied to $\mathfrak F$ and Schlichting's theorem,  there is an ideal $I$ of $D$ which is commensurable with a finite intersection $\lambda_1 H\cap \ldots \cap \lambda_n  H$ of subgroups from $\mathfrak F$. As $I$ must be trivial and $\lambda_1 H\cap \ldots \cap \lambda_n  H$ is a vector space over $C_D(a)$,  it is trivial as well. Moreover, the subvector spaces $\lambda_iH$ have codimension one in $D$ and thus $D$ has finite dimension over $C_D(a)$. As by assumption, the centralizer of $a$ has finite dimension over its center, we obtain that $D$ has finite dimension over an infinite subfield and thus over its own center  by \cite[Theorem 15.8]{lam}.
 \qed

\begin{prop}\label{Prop_RosFinCon}
A rosy division ring which has only finitely many non-central conjugacy classes is commutative.
\end{prop}

\pf
Let $D$ be a non-commutative  rosy division ring, and suppose that $D$ has only finitely many conjugacy classes.

The first step is to show that $Z(D)$ is contained in $b^D-b$ for any
non-central element $b$ in $D$. To do so, set $H$ to be  $Z(D)\cap (b^D-b)$ and observe that it is an additive subgroup of $Z(D)$. Now, following the lines of the proof of \cite[Theorem 5.6.12]{wagner} we show that $H$ has only finitely many $Z(D)$-translates. For any $z$ in $Z(D)$ we have that
$$zH = z[(b^D)-b]\cap Z(D) =  [(zb^D)-zb]\cap Z(D) .$$
Hence, if $z$ and $z'$ are two elements in $Z(D)$ such that $z'b=(zb)^d$ for some $d$ in $D$, we have that
$$
zH = (zH)^d = [(zb)^{Dd}-(zb)^d] \cap Z(D) =  [(z'b)^{D}-(z'b)] \cap Z(D) = z' H.
$$
As $Z(D)b$ contains only finitely many conjugacy classes, the group $H$ has finitely many multiplicative $Z(D)$-translates.

Now, observe that for any
two central elements $z$ and $z'$, their difference $z-z'$ belongs to $b^D-b$ if
and only if there is some element $x$ from $D^\times$ such that
$$b+z=b^x+z'=(b+z')^x.$$ As there are only finitely many conjugacy classes in $b + Z(D)$, the index of $Z(D)\cap(b^D-b)$ in $Z(D)$ has to be finite. Thus, the finite
intersection of all its $Z(D)^\times$-translates, which forms an ideal of $Z(D)$, has
finite index in $Z(D)$ as well. If $Z(D)$ is finite, the characteristic of
$D$ is positive and thus by Proposition \ref{PropCed} and Wedderburn's little theorem \cite[Theorem 13.1]{lam}
$D$ must be commutative. So we may assume that $Z(D)$ is an infinite field and hence
equal to $Z(D)\cap (b^D-b)$. Thus $Z(D)$ is contained in $b^D-b$ for any non-central element $b$ of $D$.

Now, by Kaplansky's theorem \cite[Theorem 15.15]{lam}, we can find an element $a$ in $D$ for which none of its powers belong to $Z(D)$. As $D$ satisfies the chain condition on centralizers, after replacing $a$ by one of its powers, we may assume that $C_D(a) = C_D(a^n)$ for any natural number $n$.

Now, for any natural number $m$ any element of the center has a $m$-root in $D$ by Lemma \ref{Lemma_ZDiv}. In particular, there is an infinite sequence $-1 = \xi_0, \xi_1, \xi_2, \dots$ of elements in $D$ with $\xi_{k}^{2^k} = -1$ for all $k < \omega$. It is clear that all these roots of unity have different conjugacy classes and hence all but finitely many must belong to the center since there are only finitely many non-central conjugacy classes. Furthermore, one can find two different natural numbers  $i$ and $j$ such that $\xi_i$ and $\xi_j$ belong to the center and $\xi_i a$ and $\xi_j a $ are conjugates. Thus, there is some $x$ in $D \setminus C_D(a)$ and some non-trivial $\zeta$ in the center with $\zeta^m=1$ for some $m <\omega$ such that $a=\zeta a^x$. Hence $a^m = (\zeta a^x)^m = (a^m)^x$ and so $x$ belongs to the centralizer of $a^m$ which, by the choice of $a$, coincides with the centralizer of $a$. This yields the final contradiction. \qed

The above result yields the following consequence for division rings definable in a simple theory. Recall that in a group $G$ definable in such a theory, an element $g$ of $G$ is {\em generic} over $A$ if for any $h$ in $G$ with $g$ independent from $h$ over $A$, the product $g\cdot h$ is independent from $h,A$. These elements play an essential role in the study of definable groups in this context. The next result focuses on division rings with generic elements of weight one, i.e. given any two independent elements $b$ and $c$ and any generic element $a$ either $a$ is independent from $b$ or from $c$.  A key property is that an element is generic for the multiplicative group if and only if it is generic for the additive group. For a detailed exposition, we refer the reader to \cite[Chapter 4]{wagner}.

\begin{cor}\label{CorWeight1}
A division ring with a simple theory and a generic of weight one is a field.
\end{cor}
\pf
Let $D$ be such a division ring and let $g$ be any non-central element. We denote by $\ulcorner g^D\urcorner$ the canonical parameter of the conjugacy class of $g$ in $D$. Now, let $X$ be the set of non generic elements of $D$ over $\ulcorner g^D\urcorner$. By \cite[Remark 1.1]{KruPil}, or its proof, the set of non generic elements over any given small subset forms a subdivision ring. As conjugation is an automorphism of $D$ which fixes $\ulcorner g^D\urcorner$, a such ring is invariant under conjugation and hence, we have that the division ring of non-generics over $\ulcorner g^D\urcorner$ is contained in $Z(D)$ by the Cartan-Brauer-Hua Theorem \cite[Theorem 13.17]{lam}. In fact, the division ring of non-generics over $\ulcorner g^D\urcorner$ and $Z(D)$ coincide. So $g$ itself is a generic element of $D$ independent from $\ulcorner g^D\urcorner$. Thus for any noncentral element $g$ in $D$, we have that the canonical parameter $\ulcorner g^D\urcorner$ is algebraic over the empty set. Hence $D$ has only finitely many non-central conjugacy classes and whence it is commutative by Proposition \ref{Prop_RosFinCon}.
\qed

\section{Division rings of  finite burden}

In this section we study division rings in which definable subgroups satisfiy a generalized chain condition. More precisely, given a division ring $D$ and a natural number $n$, we introduce the following property:
\begin{enumerate}
\item[$(\dagger)_n$] For any  definable subgroups $H_0,\ldots,H_n$ of $D^+$, there exists some $j\le n$ such that $\bigcap_{i\le n} H_i$ has finite index in $\bigcap_{i\neq j} H_i.$
\end{enumerate}

The motivation to analyze division rings fulfilling this property for some natural number $n$ originates in the study of division rings of finite burden, and it is as well satisfied by any superrosy division rings of finite rank. In fact, all examples of fields mentioned in the introduction satisfy $(\dagger)_1$. Below we give the precise definition of  burden.

\begin{defn} \label{Def_FinBur} Let $\pi(x)$ be a partial type.
An \emph{inp-pattern of depth $\kappa$} in $\pi(x)$ is a family of formulas $\{\psi_\alpha (\bar x;\bar y_\alpha)\}_{\alpha < \kappa}$, an array of parameters $(\bar a_{\alpha,j})_{\alpha < \kappa , j < \omega}$ with $|\bar a_{\alpha,j}|= |\bar y_\alpha|$, and a sequence of natural numbers $(k_\alpha)_{\alpha < \kappa}$ such that:
\begin{itemize}
\item the set $\{ \psi_\alpha(\bar x,\bar a_{\alpha, j})\}_{j< \omega}$ is $k_\alpha$-inconsistent for every $  \alpha < \kappa$;
\item the set $\pi(x) \cup \{ \psi_\alpha(\bar x,\bar a_{\alpha, f(\alpha)})\}_{\alpha <\kappa}$ is consistent for every $f: \kappa \rightarrow  \omega$.
\end{itemize}
We say that a theory has \emph{burden} $n$ for some natural number $n$, if there is no inp-pattern  of depth $n$ in the  partial type $x=x$. A theory of burden $1$ is called \emph{inp-minimal}. Accordingly,  a definable division ring has burden $n$ if the formula which defines the group or division ring seen as a partial type has burden at most $n$.
\end{defn}

The following result is an immediate consequence of \cite[Proposition 4.5]{cks} in the definable context.

\begin{fact}\label{Lem_Hbdn}
A definable division ring of burden $n$ satisfies $(\dagger)_n$.
\end{fact}

\begin{lemma}\label{LemmaBDNn}
A division ring satisfying $(\dagger)_n$ has dimension at most $n$ over any infinite definable subfield.
\end{lemma}
\pf Let $D$ be a division ring satisfying $(\dagger)_n$ with an infinite definable subfield $K$, and assume that the dimension of $D$ over $K$ is at least $n+1$. Choose $K$-linearly independent elements $e_0,\ldots,e_n$ in $D$. For $j\le n$, consider the definable $K$-vector spaces $V_j=\bigoplus_{i\ne j} Ke_i$. Therefore, the condition $(\dagger)_n$  yields the existence of some $k\le n$ such that the index
$$\left[ \bigcap_{j\neq k} V_j : \bigcap_{j} V_j \right] = \left[Ke_k : \{ 0 \} \right]$$
is finite, a contradiction. \qed

\begin{cor}\label{Cor_bdnn}
Any infinite division ring satisfying $(\dagger)_n$ has dimension at most $n$ over its center.
\end{cor}

\pf
Let $D$ be a division ring satisfying $(\dagger)_n$. 
By Kaplansky's theorem \cite[Theorem 15.15]{lam}, we may assume that $D$ has an element $d$ of infinite order.
Hence $Z(C(d))$ is an infinite definable subfield of $D$, so Lemma \ref{LemmaBDNn} implies that $D$ has finite dimension over $Z(C(d))$ and whence it has finite dimension over its center by \cite[Theorem 15.8]{lam}. Therefore $Z(D)$ must be infinite.
Now, we can apply Lemma \ref{LemmaBDNn} to the center of $D$ and obtain the desired result.
\qed

Immediately we obtain:

\begin{cor}\label{CorBDN}
A definable division ring of burden $n$ has dimension at most $n$ over its center. In particular, an inp-minimal division ring is commutative.
\end{cor}

Moreover, as the quaternions are a finite extension of the inp-minimal field $\mathbb R$, they have finite burden. As they are non-commutative, one cannot expect to improve the above results to obtain commutativity.

Another consequence of these results is the descending chain condition among definable subfields.

\begin{cor}\label{dcc}
Let $D$ be an infinite division ring satisfying $(\dagger)_n$. Then any descending chain of definable infinite subfields has length at most $\lfloor \log_2(n) \rfloor+1$.
Therefore, if $\F$ is a family of definable subfields of $D$, the intersection of all subfields in $\F$ is equal to a finite subintersection and so, it is definable.
\end{cor}

Now, we aim to study definable groups of automorphisms of fields satisfying $(\dagger)_n$ for some $n$. In particular, this applies to fields of finite burden as well as supperrosy fields of finite rank. We obtain results of the same nature to the following one of Hrushovski in the superstable case \cite[Proposition 3]{Hru}:

\begin{fact}\label{Hru_SSField} Any definable group of automorphisms acting definably on a definable superstable field is trivial.
\end{fact}

\begin{prop}\label{thmfinbdnGfin}
If $F$ is a field satisfying $(\dagger)_n$ and the algebraic closure of the prime field of $F$ in $F$ is infinite, then any definable group of automorphisms acting definably on $F$ has size at most $n$.
\end{prop}
\pf  Assume, as we may, that our structure is sufficiently saturated. Let $G$ be a definable group of automorphisms of $F$ acting definably, let $k$ be the prime field of $F$ and let $G_x$ denote  the stabilizer of any element $x\in k^{\rm alg} \cap F$ in $G$.
As $k$ is fixed by the action of every element in $G$ and $G/G_x$ is in bijection with the orbit of $x$, the stabilizer $G_x$ has finite index in $G$.
Now, we work with the subgroup  $$H=\bigcap_{x\in k^{\rm alg}\cap F} G_x$$ of $G$.
Note that it is a type-definable subgroup of $G$ of bounded index.
We consider the intersection $\Fix(H) = \bigcap_{\sigma \in H} \Fix(\sigma)$ of definable subfields of $F$.
By Corollary \ref{dcc} it is equal to a finite subintersection. Hence, as additionally $\Fix(H)$ contains the infinite field $k^{\rm alg}\cap F$, it is a definable infinite subfield of $F$.
Thus, Lemma \ref{LemmaBDNn} yields that $F$ has at most dimension $n$ over $\Fix(H)$, so $H$ is finite. Hence the group $G$ is a bounded definable group and whence finite by compactness. Now, consider the definable field $\Fix(G)$. By Galois theory we know that $F$ is a finite field extension of $\Fix(G)$ of degree $|G|$, and so $\Fix(G)$ is an infinite definable subfield of $F$. Hence $F$ has dimension at most $n$ over $\Fix(G)$ by Lemma \ref{LemmaBDNn} and whence $G$ has size at most $n$.
\qed

\begin{cor}
If $F$ is a field satisfying $(\dagger)_1$ and the algebraic closure of the prime field of $F$ in $F$ is infinite, then any definable group acting definably on $F$ as automorphisms is trivial.
\end{cor}

Observe that if $F$ is Artin-Schreier closed, then the algebraic closure of the prime field of $F$ is infinite in $F$. Thus, the above result holds for any infinite field of positive characteristic with finite burden and which in addition is NIP \cite{KSW} or even $n$-dependent \cite{Nadja}. 

We conclude the section with the following result in characteristic zero.

\begin{prop}\label{PropDefAut}
There are at most $n$ many definable automorphisms of a definable field of characteristic zero satisfying $(\dagger)_n$. 
\end{prop}
\pf Let $K$ be a field satisfying $(\dagger)_n$ and let $\mathcal H$ be the family of all definable automorphisms. For $\sigma$ in $\mathcal H$, let $F_\sigma$ denote the fixed field of $\sigma$, which is definable. By Corollary \ref{dcc}, the intersection $F$ of all these fixed fields is again definable. Thus, it is infinite since the characteristic is zero and hence $K$ has  dimension at most $n$ over $F$ by Lemma \ref{LemmaBDNn}. Therefore ${\rm Aut}(K/F)$ has size at most $n$ and so does $\mathcal H$ since any of its automorphisms fixes $F$.\qed

\section{Interpretable division rings in o-minimal structures}

As pointed out in the introduction, Peterzil and Steinhorn \cite[Theorem 4.1]{PS} showed that an infinite definable ring without zero divisors in an o-minimal structure is an algebraically closed field, a real closed field or the division ring of the quaternions over a definable real closed field.

Next, we generalize the above result to interpretable division rings in o-minimal structures. Different to the proof of Peterzil and Steinhorn, which has a topological flavor, our approach uses the previous sections together with more recent results on groups in o-minimal structures.

\begin{theorem}
An infinite interpretable division ring in an o-minimal structure is an algebraically closed field, a real closed field or the division ring of the quaternions over a definable real closed field.
\end{theorem}
\pf Let $D$ be a division ring interpretable in an o-minimal structure and let $K$ be its center, which is an interpretable field. As an o-minimal structure has burden one, the structure $(D,+,\times)$ has finite burden and so $D$ is finite dimensional over $K$ by Corollary \ref{CorBDN}. Consequently, the field $K$ is  infinite and interpretable in an o-minimal structure. Alternatively, any o-minimal structure is superrosy, a class of structures which is preserved under interpretation, and thus we obtain the same conclusion by Theorem \ref{ThmMain}.

By Frobenius' theorem \cite[Theorem 13.12]{lam} it suffices to show that $K$ is either algebraically closed or real closed.  To do so, we show that $K$ is isomorphic to some definable field in an o-minimal structure. Indeed, by a result of Pillay \cite{Pil}, namely that any infinite definable  field in an o-minimal structure is real closed or algebraically closed, the former is enough to conclude.

Now, consider the semidirect product $K^+ \rtimes K^\times$. 
As this group is clearly interpretable in the given o-minimal structure, it is (definably) isomorphic to a definable group $G$ by \cite[Theorem 8.23]{EPR}. This isomorphism yields the existence of a normal subgroup $N$ of $G$ isomorphic to $K^+$ and a subgroup $H$ of $G$ isomorphic to $K^\times$ such that $G=NH$ and $N\cap H=\{1_G\}$. In other words, the group $G$ is the (inner) semidirect product of $N$ and $H$ where $H$ acts on $N$ by conjugation. As $K^+ \rtimes K^\times$ is centerless, so is $G$ and thus, since both $N$ and $H$ are abelian, we obtain that $N=C_G(N)$ and $H=C_G(H)$. By the descending chain condition on definable subgroups in groups definable in an o-minimal structure, the subgroups $N$ and $H$ are definable as a finite intersection of centralizers. Now, let $e$ be the element of $N$ corresponding to $(1,1)$ in $K^+\rtimes K^\times$ under the above mentioned isomorphism. As $K^\times$ acts regularly on $K^+\setminus\{0\}$, so does $H$ on $N\setminus\{1_G\} $. Thus, for every non-trivial element $a$ in $N$ there is a unique element $h$ in $H$ such that $h * e=a$, which yields a uniquely determined identification of $N\setminus\{0\}$ and $H$. Using this, one can define a multiplication on $N$. Namely, for non-trivial elements $a_1,a_2, a_3$ in $N$, and $h_1, h_2,h_3$ their corresponding elements of $H$, we set
$$
a_1\cdot a_2 = a_3 \ \Leftrightarrow \ h_1 h_2 =h_3.
$$
This well-defined operation, together with the given group structure of $N$ seen as addition, yields a definable field structure on $N$. Thus $N$ turns out to be an infinite definable field in an o-minimal structure, which is by construction isomorphic to $K$. This finishes the proof.\qed

\end{document}